\documentstyle[11pt]{article}
\input{amssym.def}
\input{amssym.tex}

\newcommand{\Ab}{{\cal{A}{\it b}}}
\newcommand{\A}{{\cal{A}}}
\newcommand{\Bb}{{\cal{B}}}

\newcommand{\Nat}{{\cal{N}{\!\it at}}}
\newcommand{\Fac}{{\cal{F}{\!\it ac}}}
\newcommand{\NatS}{{{\cal{N}{\!\it at}}S}}
\newcommand{\FacS}{{{\cal{F}{\!\it ac}}S}}
\newcommand{\FacM}{{{\cal{F}{\!\it ac}}M}}

\newcommand{\Hom}{{\rm Hom}}
\newcommand{\Ext}{{\rm Ext}}

\newcommand{\Z}{{\bf Z}}
\newcommand{\D}{{\bf D}}
\newcommand{\B}{{\bf B}}
\newcommand{\F}{{\bf F}}
\newcommand{\G}{{\bf G}}
\newcommand{\K}{{\bf K}}
\newcommand{\Pp}{{\bf P}}
\newcommand{\Ll}{{\bf L}}

\newcommand{\up} [2]{\stackrel{#2}{#1}}
\newcommand{\too}{\longrightarrow}
\newcommand{\toop}{\up{\too}{\mbox{\large .}}}

\newcommand{\sz}[1]{\mbox{\scriptsize ${#1}$}}
\newcommand{\Arr}{-\!\!\!-\!\!\!-\!\!\!-\!\!\!-\!\!\!\to}
\newcommand{\ARr}{-\!\!\!-\!\!\!-\!\!\!-\!\!\!-\!\!\!-\!\!\!-\!\!\!-\!\!\!-\!\!\!-\!\!\!\to}

\newcommand{\sss} [3]{{#1}_{#2},\ldots ,{#1}_{#3}}
\newcommand{\msss} [3]{{#1}_{#2}\cdots {#1}_{#3}}
\begin{document}

\title{Semigroup cohomology as a derived functor}

\author{ A.\,A.\,Kostin, B.\,V.\,Novikov}
\date{}

\maketitle
\begin{abstract}
In this work we construct an extension for the category of
$0$-modules by analogy with \cite{ref5}. The $0$-cohomology functor
becomes a derived functor in the extended category. As an
application of this construction we calculate the cohomological
dimension of so-called $0$-free monoids.
\end{abstract}

{\bf 1.} 0-cohomology of semigroups appeared  in research of
projective representations of semigroups \cite{ref1}. Besides, it
was useful in studying of matrix algebras \cite{ref3} and Brauer
monoids \cite{ref4} (see also survey \cite{ref2} and references
there).

However the further study of its properties is complicated.
One of the reasons is that the semigroup 0-cohomology is not a
derived functor in the category where it is built (so-called category
of 0-modules).

The purpose of this paper is to describe the  extension
 of 0-cohomology on a larger category where it becomes a derived functor.
Our construction is similar to Baues theory for cohomology of
small categories \cite{ref5}. Therefore we omit some proofs
replacing them by references to \cite{ref5}.

As an example of application of our construction we prove  that a
cohomological dimension of a so-called 0-free semigroup equals
one. In particular, it follows that all projective representations
of a free semigroup are linearizable.

{\bf 2.} We begin with definitions. Let $S$ be a monoid. We may
assume that $S$ has a zero element (if not, let us join
it to $S$). By analogy with \cite{ref5} {\it the category of
factorizations in} $S$ is given as follows. The objects are all
nonzero elements of $S$ and the set of morphisms ${\rm Mor}(a,b)$
consists of all triples $(\alpha,a,\beta) \ \ \ (\alpha,\beta \in
S)$ such that $\alpha a\beta = b$. We will denote
$(\alpha,a,\beta)$ by $(\alpha,\beta)$ if this cannot lead to
confusion. The composition is defined by the rule: $(\alpha',
\beta') (\alpha, \beta)= (\alpha'\alpha,\beta\beta')$; hence we
have $(\alpha,\beta)=(\alpha,1) (1,\beta)=(1,\beta)(\alpha,1).$
Denote this category by $\FacS$.

{\it A natural system on} $S$ is a functor $\D:\FacS \too \Ab.$
The category $\NatS = \Ab^\FacS$ is an Abelian category with
enough projectives and injectives \cite{ref6}. Denote the value of
$\D$ at the object $a\in {\rm Ob}\FacS$ by $\D_a$. By $\alpha_*$
and $\beta^*$ denote values of $\D$ at morphisms $(\alpha,1)$ and
$(1,\beta)$ respectively. We have
$\D(\alpha,\beta)=\alpha_*\beta^*$ for all morphisms
$(\alpha,\beta)$.

For given natural number $n$ denote by $Ner_n S$ the set of all
$n$-tuples $(a_1,\ldots,a_n), \ a_i\in S$, such that
$\msss{a}{1}{n}\ne 0$. By definition $Ner_0 S = \{1\}.$ A {\it
n-cochain} assigns to each point $a=(a_1,\ldots,a_n)$ of $Ner_n S$ an
element on $\D_{\msss{a}{1}{n}}$. The set of all $n$-cochains is
an Abelian group $C^n(S,\D)$ with respect to the pointwise
addition. Set $C^0(S,\D)=\D_1$.

The {\it coboundary} $ \delta = \delta^n:C^n(S,\D)\too
C^{n+1}(S,\D) $ is given by the formula $(n \geq 1)$
\begin{eqnarray}
&&(\delta f)(a_1,\ldots,a_{n+1})= a_1{}_* f(a_2,\ldots,a_{n+1})
\nonumber \\ &&+\sum_{i=1}^{n}{(-1)^i f(a_1,\ldots,a_i
a_{i+1},\ldots ,a_{n+1})} +(-1)^{n+1} a_{n+1}^*f(a_1,\ldots, a_n).
\nonumber
\end{eqnarray}

For $n=0$ let $\delta^0:C^0(S,\D)\too C^1(S,\D)$ be defined by
$$
\delta f(x)= x_*f-x^*f\ \ \  (f\in D_1,\ \ x\in S\setminus 0).
$$
One can check directly that $\delta^n\delta^{n-1}=0.$ By
$H^n(S,\D)$ denote the cohomology groups of the complex
$\{C^n(S,\D),\delta^n\}_{n\geq0}.$

{\bf 3.} Now we define a {\it trivial natural system} $\Z$. To each
object $a\in S\setminus 0$ it
assigns the infinite cyclic group $\Z_a$ generated by a symbol
$[a]$; and to each morphism $(\alpha,\beta):a\too b$ it assigns  a
homomorphism of the groups $\Z(\alpha,\beta):\Z_a\too\Z_b$
which takes $[a]$ to $[b]$.

Since $\NatS$ has enough projective and injective, hence there
exists the derived functor $\Ext^n_\NatS(\Z,-)$. This functor is
isomorphic to the co\-ho\-mo\-lo\-gy functor $H^n(S,-)$ which is
defined in Section {\bf 2.} To prove this statement we construct a
suitable projective resolution of $\Z$.

For every $n\geq 0$ we denote by $\B_n:\FacS\too \Ab$ the
following natural system. For an object $a\in S\setminus 0$ the
group $\B_n(a)$ is a free Abelian group generated by the set of
symbols $[\sss{a}{0}{n+1}]$ such that $\msss{a}{0}{n+1}=a.$ To
each morphism $(\alpha,\beta)$ we assign a homomorphism of groups by
the formula
$$
\B_n(\alpha,\beta):\ [\sss{a}{0}{n+1}]\longmapsto
[\alpha \sss{a}{0}{n+1}\beta].
$$
The functors $\B_n\ \ (n\geq 0$)
constitute a chain complex $\{\B_n,
\partial_n \}_{n\geq 0}$, where $\partial_n:\B_n\toop \B_{n-1}\
\ (n\geq 1)$ is a natural transformation with the set of its
components
$$
(\partial_n)_a:\B_n(a)\too \B_{n-1}(a),
$$
$$
(\partial_n)_a[\sss{a}{0}{n+1}]=
\sum_{i=0}^{n}{(-1)^i[a_0,\ldots,a_ia_{i+1},\ldots,a_{n+1}].}
$$

{\bf 4.} {\sc Lemma. }{\it The natural system $\B_n$ is a
projective object in $\NatS.$
}

\noindent
{\sc Proof. } Consider the following diagram with the exact row $$
\begin{array}{ccccc}
&&{}\,\B_n        &      &  \\
 &&\hspace{0.5em}\Biggr\downarrow\sz{\nu}& &\\
 \D\hspace{-1em} & \up{\Arr}{\mu}\hspace{-1em} & {\bf E} & \hspace{-1em}\Arr &\hspace{-1em} 0
\end{array}
$$
%%%
%%%    Diagrama (1) p.5
%%%
and construct a natural transformation $\tau:\B_n\toop\D$ which
turns this diagram into commutative.

Let $s=\msss{s}{0}{n+1},\ \ \hat{s}=\msss{s}{1}{n}.$ Choose
$a_{(\sss{s}{1}{n})}\in \D(\hat{s})$ such that
$\mu_{\hat{s}}a_{(\sss{s}{1}{n})}=\nu_{\hat{s}}[1,\sss{s}{1}{n},1]$,
and put
$$
\tau_s[\sss{s}{0}{n+1}]=\D(s_0,s_{n+1})a_{(\sss{s}{1}{n}).}
$$

The natural transformation is well defined. Indeed,
\begin{eqnarray}
&&\tau_{\alpha s\beta}\B_n(\alpha ,\beta)[\sss{s}{0}{n+1}]=
\D(\alpha s_0,s_{n+1}\beta)a_{(\sss{s}{1}{n})}=\nonumber \\
&&\D(\alpha,\beta)\D(s_0,s_{n+1})a_{(\sss{s}{1}{n})}=
\D(\alpha,\beta)\tau_s[\sss{s}{0}{n+1}].\ \ \Box\nonumber
\end{eqnarray}

{\bf 5.} {\sc Lemma. }{\it The chain complex
$\{\B_n,\partial_n\}_{n\geq 0}$ is a projective resolution of the
natural system $\Z$. }

The proof is similar to \cite{ref5}.

{\bf 6.} Now we are ready to prove the main result of this
paper.\\ {\sc Theorem. }{\it For any monoid $S$ with a zero
element  there is an isomorphism of the functors:
$$
H^n(S,-)\cong\Ext_{\NatS}^n(\Z,-).
$$ }
\\
{\sc Proof. } Define an isomorphism of complexes
$$
\Psi_{\D}^*:\{\Hom_{\NatS}(\B_n,\D),\partial^n\}_{n\geq 0}\too
\{C^n(S,\D),\delta^n\}_{n\geq 0}
$$
(here we denote
$\partial^n=\Hom_{\NatS}(\partial_{n-1},\D)$) as follows. Let the
homomorphism of Abelian group
$$
\Psi_{\D}^n:\Hom_{\NatS}(\B_n,\D)\too C^n(S,\D)
$$
be given by
$$
(\Psi_{\D}^n\tau)(\sss{a}{1}{n})=
\tau_{\msss{a}{1}{n}}[1,\sss{a}{1}{n},1]\in \D_{\msss{a}{1}{n}} \
{\rm for}\ \msss{a}{1}{n}\ne 0.
$$
Let $a=\msss{a}{0}{n+1}$, i.e.
$[\sss{a}{0}{n+1}]\in \B_n(a).$ Since the diagram
$$
\begin{array}{ccc}
\B_n(\msss{a}{1}{n})&\hspace{-2ex}\up{\Arr}{\tau_{\msss{a}{1}{n}}}&\hspace{-2ex}\D_n(\msss{a}{1}{n})\\
\hspace{-10ex}\sz{\B_n(a_0,a_{n+1})}\Biggr\downarrow&&\Biggr\downarrow\sz{\D_n(a_0,a_{n+1})}\hspace{-7ex}\\
 \B_n(a) &\hspace{-1.5em} \up{\ARr}{\tau_a}&\hspace{-1.5em}\D_n(a)
\end{array}
$$
%%
%%Diagram on p.6
%%
is commutative we have
$$ \tau_a[\sss{a}{0}{n+1}]=
\D(a_0,a_{n+1})\tau_{\msss{a}{1}{n}}[1,\sss{a}{1}{n},1].
$$
Therefore $\Psi^n_{\D}\tau=0$ implies that $\tau_a$ vanishes on
all generators of the group $\B_n(a).$ Hence $\Psi_{\D}^n$ is
injective.

Further, for any $f\in C^n(S,\D)$ define a natural transformation
$\varphi:\B_n\toop \D$:
$$
\varphi_a[\sss{a}{0}{n+1}]=
\D(a_0,a_{n+1})f(\sss{a}{1}{n})
$$
It is clear that
$\Psi^n_{\D}\varphi=f$ and hence $\Psi^n$ is surjective. The
commutativity of the diagram  $$
\begin{array}{ccc}
\Hom_{\NatS}(\B_n,\D)&\hspace{-2ex}\up{\Arr}{\partial^n}&\hspace{-2ex}\Hom_{\NatS}(\B_{n+1},\D)\\
\hspace{-3ex}\sz{\Psi^n_{\D}}\Biggr\downarrow&&\Biggr\downarrow\sz{\Psi^{n+1}_{\D}}\hspace{-0ex}\\
 C^n(S,\D)&\hspace{-1.5em} \up{\ARr}{\delta^n}&\hspace{-1.5em}C^{n+1}(S,\D)
\end{array}
$$
%%
%% The diagram on p.7
%%
is established immediately.

It can easily be checked that the family $\Psi^n=\{\Psi_{\D}^n |
\D\in\NatS\}$ is a natural transformation. From above we see that
$\Psi^n$ induces an isomorphism of functors $H^n$ and $\Ext^n$. $
\Box$

{\bf 7.} Let us discuss the relation between cohomology which is
defined above and cohomology groups of other kinds. In Section
{\bf 1} we note that the 0-cohomology is a particular case of our
construction. This can be shown in the following way. Let $A$ be
an Abelian group and  $\bf A$ be a natural system  given by
$$
{\bf A}(s)=A\ \ {\rm and} \ \ \alpha_*\beta^*a=\alpha a
$$
for all $s\in \FacS, \ \ (\alpha,\beta)\in {\rm Mor}\FacS.$ In other
words, {\bf A} is so-called 0-module over S \cite{ref1}: an action
$(S\setminus\{0\})\times A\too A$ is given, which satisfies the
following conditions:
$$
s(a+b)=sa+sb,
$$
$$ st\ne 0 \Rightarrow s(ta)=(st)a,
$$
where $s,t\in S\setminus 0$ and $a,b\in A.$ 0-Cohomology groups
are denoted by $H_0^n(S,A)$.

Note that Eilenberg-MacLane cohomology of semigroups \cite{ref8}
can be considered as a particular case of the 0-cohomology. Namely, if
$S$ is a semigroup (possibly without a zero), then $H^n(S,-) \cong
H_0^n(S^0,-)$, where $S^0$ is the semigroup $S$ with an adjoint
zero.

The category of 0-modules arises naturally in  applications of
0-co\-ho\-mo\-logy theory \cite{ref4}. However it is easily
shown that the second 0-cohomology group of the commutative
semigroup $S=\{u,v,w,0\}$ with $u^2=v^2=uv=w,\ \ uw=vw=0$ (see
\cite{ref1}) is not trivial for all nonzero 0-module over $S$.
Therefore the 0-cohomology is not a derived functor on the category of
0-modules. This is the reason for introducing the category $\NatS$.

Our construction differs from Baues' cohomology theory for monoids
\cite{ref5} in the  first step only. Actually in \cite{ref5} a
monoid $S$ is regarded as a category with a single object. At the
same time the Baues' category of factorizations in $S$ is equal to
$\Fac S^0$ out of Section {\bf 2}. Therefore the Baues' cohomology groups
of monoid $S$ and cohomology
grops of $S^0$ in our sense are the same. However if
$S$ possesses a zero element then the category $\FacS$ and Baues'
one are not equivalent and we obtain the different cohomology
groups.

{\bf 8.} Let us consider an application of the obtained results.
{\it Cohomological dimension} {\rm c.d.}$S$ of monoid $S$ is
the greatest natural number such that $H^n(S,\D)\ne 0$ for some
$\D\in \NatS.$ The Theorem from Section {\bf 6} allows us to use a
projective resolution for calculation of the dimension.

It is well-known that in many cohomological theories  c.d. of
free objects equals 1. Free objects in the class of monoids with
zero are free monoids with adjoint zero element. Nevertheless in
our case the family of monoids having c.d.1 is larger.

A monoid is called a {\it 0-free monoid} if it is isomorphic to a
Rees factor monoid of a free monoid. Free monoids with adjoint
zero will be regarded as 0-free monoids too.

{\bf 9. } We shall need the following\\ {\sc Lemma. }{\it Let
$\A,\Bb$ be categories, $\F:\A\too \Bb$,
$\G:\Bb\too\A$ be adjoint functors $(\F\dashv\G),$  functor $\G$
preserves epimorphisms and the counit $\varepsilon:\F\G\toop {\rm Id}_{\Bb}$
is identical. If an object
$a\in\A$ is projective then $\F(a)$ is  projective too. }\\ {\sc
Proof. } Let $a\in\A$ be a projective object. Consider a diagram
%$$
%\begin{array}{rrcll}
%   &                & {}\,\F(a)        &      &  \\
%   &                &\hspace{2ex}\rV{\alpha} &      &  \\
%c & \up{\too}{\beta} & b      & \too & 0
%\end{array}
%$$
$$
\begin{array}{rrcll}
   &                & {}\ \F(a)        &      &  \\
 &&\hspace{1ex}\Biggr\downarrow\sz{\alpha}& &\\
c\hspace{-0em} & \up{\Arr}{\beta}\hspace{-1em} & b& &
\end{array}
$$

%%
%% First diagram on p.9
%%
with the exact row $(c,b\in\Bb).$  Since functor $\G$ preserves
epimorphisms we obtain the diagram:
\begin{equation}\label{f2}
\begin{array}{rrcll}
   &                & \hspace{-1ex}a        &      &  \\
   &                &\hspace{4.4ex}\Biggr\downarrow\sz{\G(\alpha)\eta_a} &      &  \\
\G(c) & \up{\Arr}{\G(\beta)}\hspace{-2em} &\G(b) & &
\end{array}
\end{equation}
%%
%% Third  diagram on p.9
%%
where $\eta:{\rm Id}_{\A}\toop\G\F$ is the unit of the adjunction
$\F\dashv\G$. Since $a$ is projective, there is a homomorphism
$\gamma:a\too\G(c)$ which makes  diagram (\ref{f2}) commutative.
This means that $\G(\beta)\gamma =\G(\alpha)\eta_a$ and
$\beta\F\gamma =\alpha\F(\eta_a).$ Using the equalities
$\F(\eta_a)= {\rm Id}_{\F(a)}$ and $\F\G={\rm Id}_{\Bb}$ we get
$\beta\F\gamma= \alpha.\ \Box$

{\bf 10.} {\sc Theorem. }{\it $c.d.M\leq 1$ for all 0-free monoids
$M$. }\\ {\sc Proof. } For a given monoid $M$ consider the exact
sequence $$ 0\too\Pp_M\toop\B_M\toop\Z_M\too 0 $$ where $\Z_M,
\B_M$ are natural systems defined in Section {\bf 3}, $\Pp_M={\rm
Ker}(\B_M\toop\Z_M)$. We need to prove that $\Pp_M$ is a
projective functor.

It follows from Section 7 that $\Pp_M$ is a free functor whenever
$M$ is a free monoid with adjoint zero (see \cite{ref5}, Lemma
6.7).

Now let $M$ be a 0-free monoid, $M\cong W/I$ where $W$ is a free
monoid  and $I$ is an ideal in $W$. Consider the category of
factorizations $\F W$ which was defined in \cite{ref5}, i.e.
$\F W=\Fac (W^0)$. Define the functor $\K:\FacM\too \F W$ which takes
each nonzero element from $M$ to its preimage under the canonic
homomorphism $W\too W/I$. Functor $\K$ is well defined and induces
the functor $\K^*:{\bf Nat}W\too \Nat M$, where ${\bf Nat}W=\Ab^{\F W}$.

Consider the exact sequence which is defined in \cite{ref5},
Sec.5: $$
0\too\tilde{\Pp}_W\up{\too}{\tilde{\delta}_W}\tilde{\B}_W\up{\too}{\tilde{\varepsilon}_W}\tilde{\Z}_W\too0,
$$ where $\tilde{\Pp}_W, \tilde{\B}_W, \tilde{\Z}_W:\F W\too\Ab$
are natural systems on $W$. We have $$ \K^*(\tilde{\Z}_W)=\Z_M,\
\K^*(\tilde{\B}_W)=\B_M,\
\K^*(\tilde{\varepsilon}_W)=\varepsilon_M $$ hence
$\K^*(\tilde{\Pp}_W)=\Pp_M$.

Consider the functor $\Ll:\Nat M\too {\bf Nat}W$ which is given by
$$
\Ll(\G)_a=\left\{
\begin{array}{rl}
\G_a,   & {\rm if}\  a\not\in I\\ 0,
      & {\rm if}\  a\in I
\end{array}\right.
$$
where $\G\in \Nat M$, and
$$
\Ll(\G)(x,a,y)=\left\{
\begin{array}{rl}
\G(x,a,y),   & {\rm if}\  xay\not\in I\\
 0,      & {\rm if}\  xay\in I
\end{array}\right.
$$
Evidently $\K^*\Ll={\rm Id}_{\Nat M}$ and there is a natural
transformation ${\rm Id}_{{\bf Nat}W}\toop\Ll\K^*$. It implies
that $\Ll$ is right adjoint to $\K^*$. Besides, $\Ll$ preserves
epimorphisms and by \cite{ref5} $\tilde{\Pp}_W$ is a free object.
Using Lemma 9 we get $\Pp_M$ is a projective object. $\Box$

{\bf 11.} The semigroup is called {\it  0-cancellative } if $$
ax=bx\ne 0 \Rightarrow a=b \ {\rm and}\ xa=xb\ne 0 \Rightarrow a=b
$$ for all  elements $a,b,x$. In view of Theorem 10 the following
question arises: is a 0-cancellative monoid of cohomological
dimension one a 0-free monoid?

{\it E-mails:
\begin{quote} {
  andreykostin@mail.com\\
  boris.v.novikov@univer.kharkov.ua
}
 \end{quote}
}

\end{document}